\documentclass[11pt,reqno]{amsart}

\usepackage{amsmath,amsthm,amssymb, mathrsfs, bbm}
\usepackage[breaklinks=true]{hyperref}
\usepackage{color,soul} %% for editing only - remove at the end
\usepackage{graphicx}
\usepackage[utf8]{inputenc}
\usepackage[margin=1.5in]{geometry}

\theoremstyle{plain}

\theoremstyle{definition}

\usepackage[backend = biber,isbn=false,doi=true,eprint=false,url=false,style=alphabetic,maxbibnames=10,maxalphanames=4]{biblatex}

\DeclareSourcemap{
  \maps[datatype=bibtex]{
    \map{
      \step[fieldsource=author,
            match=Eschmeier,
            final]
      \step[fieldset=keywords, fieldvalue=JEarticle]
    }
  }
}
\numberwithin{equation}{section}

\newcommand{\C}{\mathbb{C}}

\begin{document}
\title{J\"org Eschmeier's mathematical work}

\author{Ernst Albrecht}
\address{Fachrichtung Mathematik, Universit\"at des Saarlandes, 66123 Saarbr\"ucken, Germany}
\email{ernstalb@math.uni-sb.de}
\author{Ra\'ul E. Curto}
\address{Department of Mathematics, University of Iowa, Iowa City, Iowa, USA}
\email{raul-curto@uiowa.edu}
\author{Michael Hartz}
\address{Fachrichtung Mathematik, Universit\"at des Saarlandes, 66123 Saarbr\"ucken, Germany}
\email{hartz@math.uni-sb.de}
\thanks{M.H. was partially supported by a GIF grant and by the Emmy Noether Program of the German Research Foundation (DFG Grant 466012782).}
\author{Mihai Putinar}
\address{University of California at Santa Barbara, CA,
USA and Newcastle University, Newcastle upon Tyne, UK} 

\email{\tt mputinar@math.ucsb.edu, mihai.putinar@ncl.ac.uk}
\thanks{M.P. was partially supported by a Simons Foundation collaboration grant.}

\begin{abstract} An outline of J\"org Eschmeier's main mathematical contributions is organized both on a historical perspective, as well as on a few distinct topics.
The reader can grasp from our essay the dynamics of spectral theory of commutative tuples of linear operators during the last half century. Some clear directions of future research are also underlined. 
\end{abstract}

\maketitle

J\"org Eschmeier grew up mathematically in the vibrant atmosphere of late XX-th Century German Modern Analysis.
His doctoral advisor, Heinz G\"unther Tillmann was a scientific grandson of Otto Toeplitz. He instilled with high competence in J\"org a lifelong
fascination with distribution theory and general duality theory in locally convex spaces.
  This happened around the late 1970-ies at the University of M\"unster,
where also George Maltese offered J\"org the challenge of an active research group.
At the same time and in the same place, function theory of several complex variables was blooming, turning M\"unster into one of the leading world centers on
the subject. J\"org was exposed early on during his studies to analytic techniques, homological algebra methods and geometric interpretations of the intricate nature of holomorphic functions of several complex variables. Throughout his brilliant career, he masterfully combined these two main streams of his student years.
We collect below a few pointers to J\"org Eschmeier's highly original discoveries. Our text merely touches upon his deep impact on contemporary operator theory, conveying a fraction of his superb scientific orientation and elegant style of pursuing research.

\section{Duality and spectral decompositions} The reader of this essay has unquestionably been fascinated by some sort of spectral decompositions. It is said that John von Neumann was asked by reporters on his death bed what the highest scientific discovery of his outstanding career was. To the stupefaction of all he mentioned the spectral theorem for unbounded self-adjoint operators, an esoteric concept for the layman, leaving aside his contributions to mathematical economics, game theory, nuclear weapons, electronic computers, fluid mechanics and much more.
\begin{figure}
\includegraphics[width=9cm]{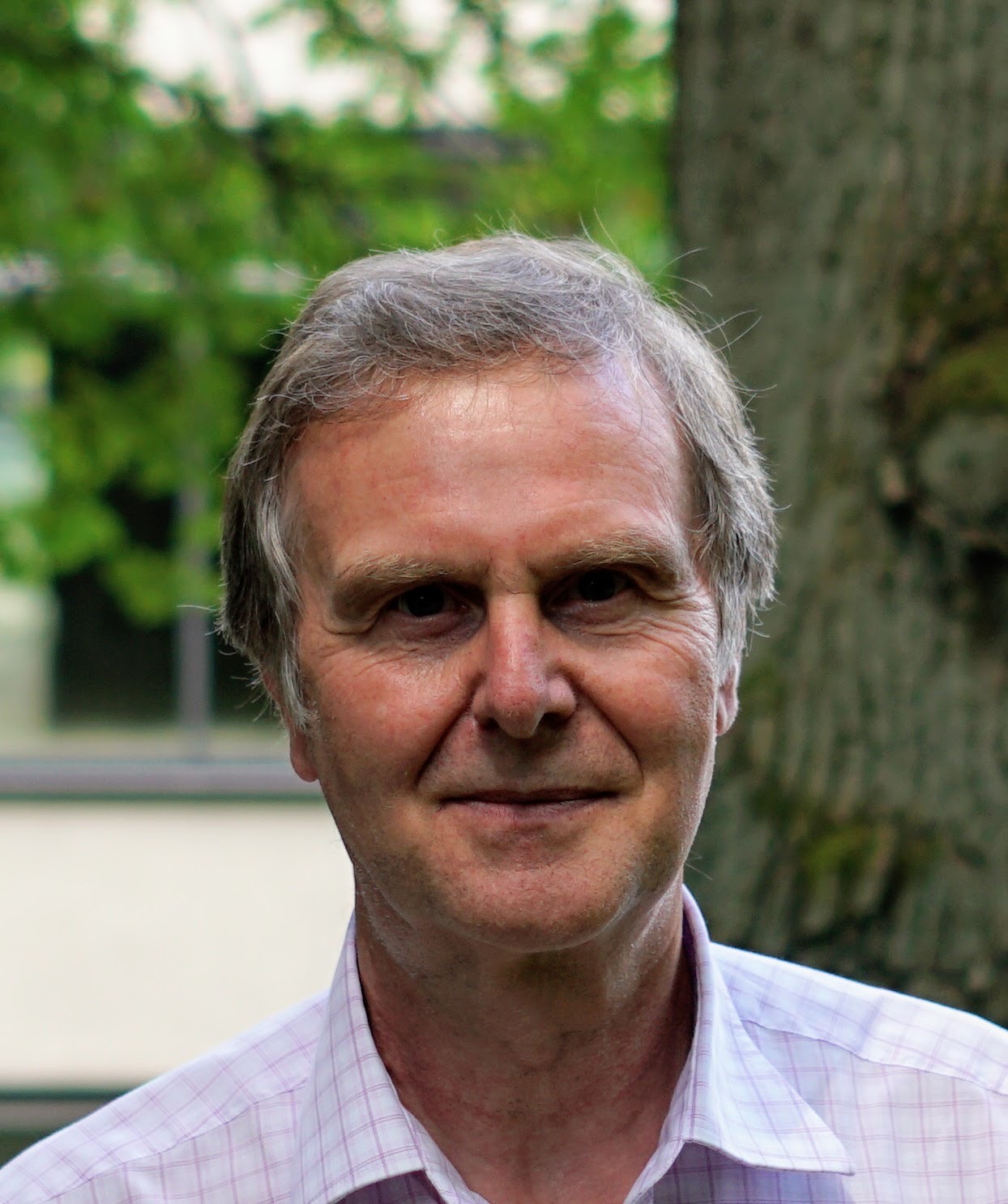}
\caption{J\"org in 2018}
\end{figure}
The ubiquitous eigenvector of a finite, self-adjoint matrix has to be replaced when dealing with infinite self-adjoint matrices carrying a continuous spectrum by a subspace of  vectors representing a localized window in the spectrum. This necessary step in mathematical spectral analysis stirred countless discussions among quantum physicists. Generalized eigenfunctions were proposed by them as a natural substitute, obtained however with the price of stepping outside the original Hilbert space. Dirac's generalized function $\delta$ stands aside in this respect.

On a totally perpendicular direction and in full resonance with the Bourbaki tendencies of the day, groups of mathematicians explored in the second part of XX-th Century several comprehensive axiomatic approaches to spectral decomposition. The third volume of Dunford and Schwartz monumental monograph \cite{Dunford-Schwartz-3} contains ample references on this forgotten chapter of operator theory. We owe to the genius of Erret Bishop \cite{Bishop} the leap forward and to the unmatched insight of Foia\c{s} \cite{Foias-1963} the foundational concept of decomposable operator.

Let $X$ be a Banach space over the complex numbers and let $T \in L(X)$ be a linear bounded operator acting on $X$. The operator $T$ is called {\it decomposable}
if, for every finite open cover of its spectrum,
$$ \sigma(T) = \sigma(T,X) \subset U_1 \cup U_2 \cup \ldots \cup U_n$$ there are 
$T$-invariant subspaces $X_1, X_2, \ldots, X_n$ with localized spectra of the respective restrictions of $T$:
$$ \sigma(T, X_j) \subset U_j, \ \ 1 \leq j \leq n,$$
and spanning the whole space: $X = X_1 + X_2 + \ldots + X_n$. Normal operators on Hilbert space, compact operators, and representations of functions algebras carrying a partition of unity are all decomposable.

%%%%%%%%%%
% Personal remarks of EA:
\begin{figure}
\includegraphics[width=12cm]{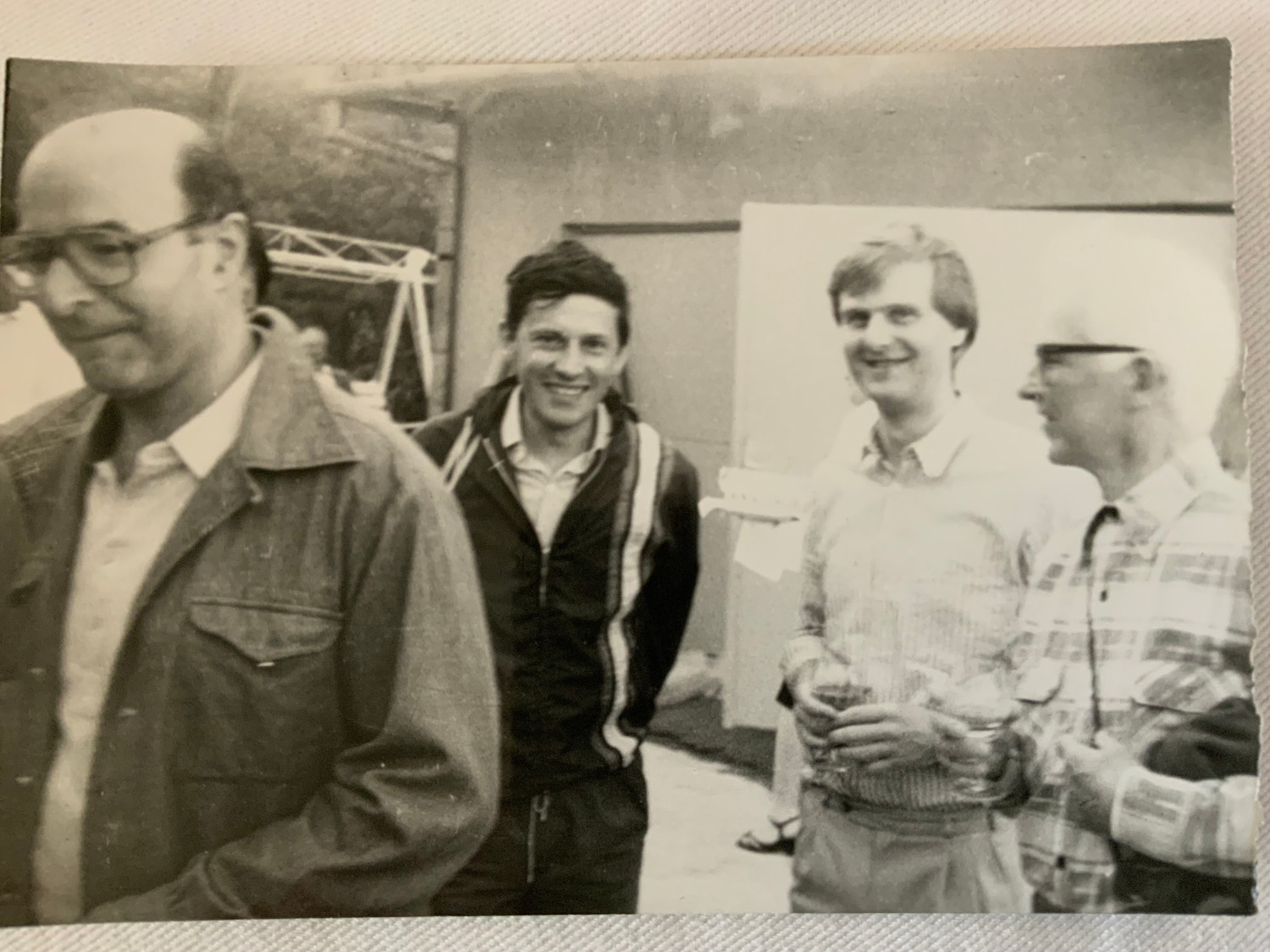}
\caption{1988 Timisoara Operator Theory Conference. From left to right: Sasha Helemskii, Mihai Putinar, J\"org Eschmeier, Henry Helson}
\end{figure}
At the beginning of a theory, variants of the definitions and permanence problems
usually have to be clarified. In the original definition of a decomposable operator Foia\c{s} in \cite{Foias-1963} required the
subspaces $X_j$ to be spectral maximal and in the last chapter of 
their monograph \cite{C-F-1968}
Colojoar\u{a} and Foia\c{s} ask if the restrictions of decomposable operators
are again decomposable. As J\"org told one of us (E.A.), he first came in touch with the
theory of decomposable operators in a student seminar of Tillman in 1977/78,
which was announced by his assistant Erich Marschall as follows:
``Here is a {\it Lecture Notes} volume \cite{E-L-1977} with a positive answer to
that question and a manuscript \cite{Albrecht-1978} with a negative answer.
Let us find out which one is correct.'' 
In \cite{JE-1988-a} J\"org published a
further class of examples which are much more natural and less technical than
the example in \cite{Albrecht-1978}. After this seminar he wrote his
(unpublished) Diploma Thesis \cite{JE-Diplom} on local decomposability
(in the sense of Vasilescu \cite{Vas-1969, Vas-1971}) and
functional calculi for closed linear operators on Banach spaces.
In particular (using methods from \cite{Albrecht-1979})
he already obtained most of the characterization of such operators described in Theorem IV.4.26 of \cite{Vas-1982}.
%%%%%%%%%%

With the proper definition of joint spectrum $\sigma(\tau,X)$ of a tuple of commuting linear bounded operator $\tau$ proposed by Taylor \cite{Taylor-1970}, the notion of decomposability carries over, but the challenging questions and complications multiply.

J\"org Eschmeier's dissertation of 1981 \cite{JE-1981}, summarized in his {\it Inventiones} article \cite{JE-1982}, is devoted to the study of spectral localization within the novel, at that time, Taylor analytic functional calculus. His lucid insight builds on solid ground well cultivated in Germany at that time, notably an earlier contribution of one of us (E.A.) \cite{Albrecht-1974}, who is incidentally a mathematical nephew of J\"org via the filiation Tillmann-Gramsch. J\"org navigates there with ease and efficiency through the formidable Cauchy-Weil integral representation formulas invoked by Taylor. He would then continue for two good decades to explore the subject.

But progress was not possible without Bishop's ideas. In a nutshell, Bishop claimed that it is not generalized eigenvectors \`a la Dirac's distribution that illuminates general spectral decomposition behavior, but analytic functionals. Fortunately, the necessary passage from Schwartz distributions to analytic functionals or hyperfunctions was much explored by mathematical analysts of the 1980-s decade. To be more specific, in order to better understand the spectral characteristics of the linear bounded operator $T \in L(X)$, Bishop points to the space of $X$-valued analytic functions defined on an open set $U$, ${\mathcal O}(U)^X$, and the linear pencil:
$$ z I - T : {\mathcal O}(U)^X \longrightarrow {\mathcal O}(U)^X.$$
The resolvent of $T$, restricted to values of $z$ outside the spectrum, is lurking around. The operator $T$ is said to have the single valued extension property
if the map $z I - T$ above is injective for every open set $U$. In that case a localized spectrum of $T$ with respect to a vector makes sense. The operator $T$ satisfies Bishop's condition $(\beta)$ if the map $z I - T$ is injective with closed range for every open set $U$. A timely 1983 observation of one of us (M.P.) \cite{Putinar-1983} interprets these notions in terms of sheaf theory:
$$ {\mathcal F}_T(U) = {\rm coker} (z I - T : {\mathcal O}(U)^X \longrightarrow {\mathcal O}(U)^X) = {\mathcal O}(U)^X)/(z I - T){\mathcal O}(U)^X,$$
is an analytic sheaf if $T$ has the single valued extension property and it is an analytic sheaf of Fr\'echet spaces if $T$ satisfies property $(\beta)$. Note that the canonical identification 
$$ {\mathcal F}_T(\C) = X, \ \ \ T = M_z,$$
takes place, with the operator $T$ represented by multiplication by the complex variable. Moreover, the spectrum of $T$ is in this case equal to the support of the sheaf
model $ {\mathcal F}_T$, the support of a section (a.k.a. vector) $x \in X$ is the local spectrum, and so on. The crucial insight came from the proof that $T$ is decomposable if and only if the canonical sheaf model (and as a matter of fact any sheaf model) is soft \cite{Putinar-1983}. Now, analytic sheaves were at home both at M\"unster and Bucharest. A fruitful collaboration between J\"org and Mihai was built on this rich ground, leading to a dozen joint publications. The rest is technique, sheaf cohomology back and forth, in the treacherous topological homological setting.

Returning to Bishop \cite{Bishop}, his clear message was to consider in parallel the non-pointwise behavior of the linear pencils $z I - T$ and $z I - T'$, where $T'$ is the topological dual of $T$. This is for the simple reason that elements of the dual of ${\rm coker} (z I - T)$ on ${\mathcal O}(U)^X$ are analytic functionals $\phi \in {\mathcal O}(U)'^{X'}$ which fulfill the generalized eigenvector equation $(z I - T')\phi = 0$. Playing this distribution theory game at the abstract multidimensional level bare fruit:
\begin{enumerate}
\item{A commutative tuple $\tau$ of linear bounded operators satisfies Bishop's property $(\beta)$ if and only if admits a quasi-coherent analytic, Fr\'echet sheaf model;}
\item{A commutative tuple $\tau$ of linear bounded operators is decomposable if and only if both $\tau$ and $\tau'$ satisfy property $(\beta)$;}
\item{A commutative tuple $\tau$ of linear bounded operators is the restriction of a decomposable tuple to a joint invariant subspace if and only if $\tau$ satisfies property $(\beta)$;}
\item{Two quasi-similar tuples of operators subject to property $(\beta)$ have equal joint spectra and equal essential joint spectra;}
\item{Division of distributions by complex analytic functions follows from a Bishop's type property with respect to smooth functions;}
\item{The functoriality of the sheaf model of a commutative tuple of linear bounded operators implies Riemann-Roch Theorem on singular analytic spaces.}
\end{enumerate}
Some of the definitive results enumerated above were obtained in several increments, each a source of joy and bewilderment \cite{JE-MP-1984,Putinar-1986,JE-MP-1988,
  JE-1992,EA-JE-1997}. The monograph \cite{JE-MP-book} collects such advances obtained until early 1990-ies and unified by the concept of analytic sheaf model. A decade later, the monograph by Laursen and Neumann \cite{LN-book} offered a complementary, cohomology free account of local spectral theory 
(in the case of single operators). 
\begin{figure}
\includegraphics[width=12cm]{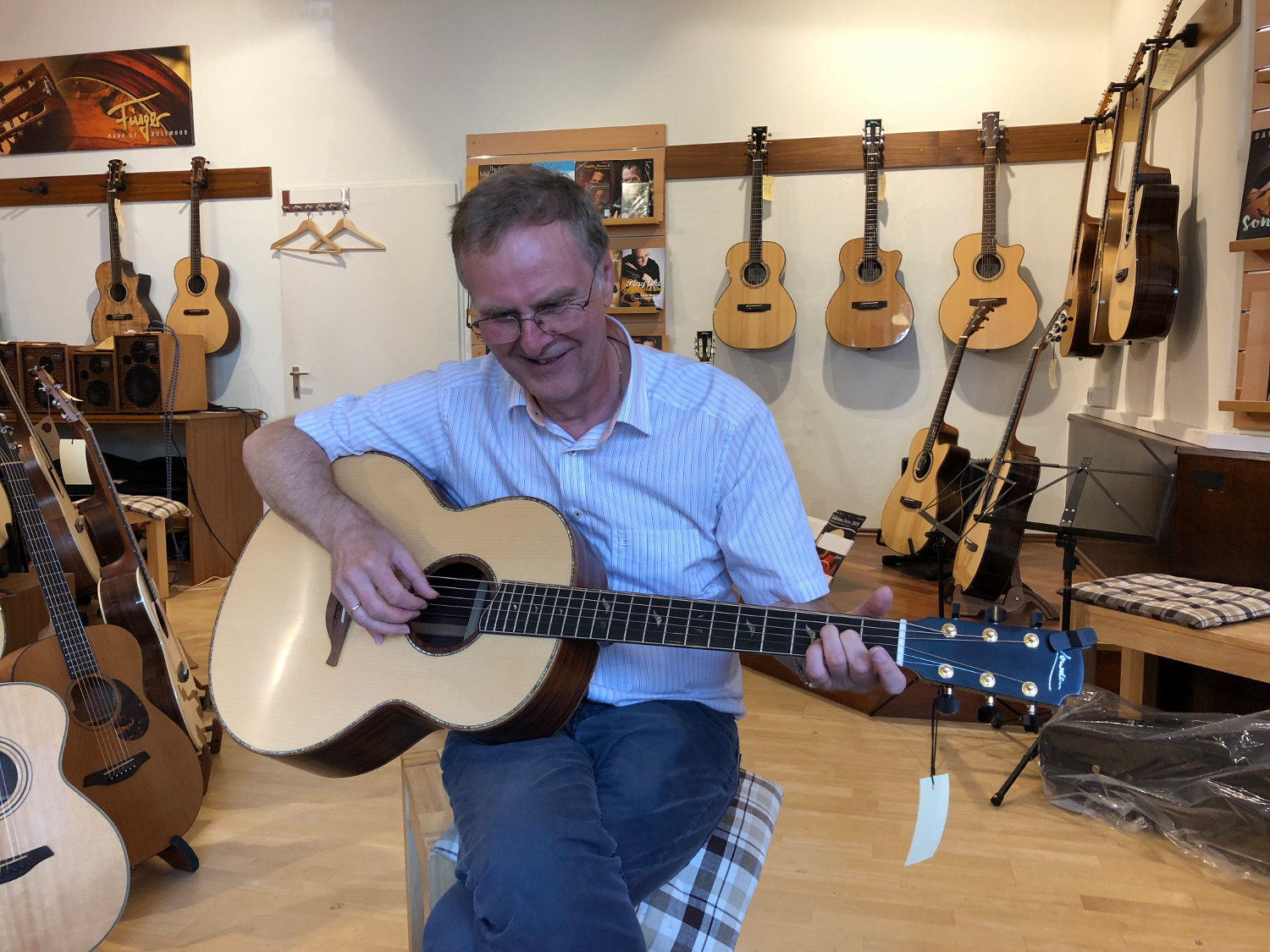}
\caption{J\"org in 2019}
\end{figure}

\section{Invariant subspaces}

The famous results and in particular the methods of proof in Scott Brown's
work \cite{Brown78,Brown87} on invariant subspaces for subnormal operators
and hyponormal operators with thick spectrum on Hilbert spaces had a great
impact for the further development in the research on the invariant subspace
problem. In \cite{Brown87}, Scott Brown used the fact (due to  \cite{P84})
that a hyponormal operator $T$ is subscalar (and hence subdecomposable) 
to obtain nontrivial invariant subspaces for $T$ when the
spectrum $\sigma(T)$ is thick in the sense that for some non empty set
$G\subset \mathbb{C}$ the set $\sigma(T)\cap G$ is dominating in $G$.
He also noticed that subdecomposability is sufficient for his result.
In \cite{JE89a} J\"org obtained a first result for general Banach spaces:
If $T\in L(X)$ is subscalar with $\text{int}({\sigma(T)})\neq \emptyset$,
then $T$ has nontrivial invariant subspaces. 

By \cite{EA-JE-1997}, subdecomposability is equivalent to Bishop's property
$(\beta)$, which is a local property and has a nice duality theory.
Using these results J\"org removed restrictions on the operators or on
the underlying Banach spaces to obtain some first localized invariant subspace
results. Finally, in joint work with Bebe Prunaru \cite{JE-Pru90}, he proved that every continuous linear operator $T$ on an arbitrary Banach
space $X\neq\{0\}$ which, for some compact set $S\subset \mathbb{C}$,
satisfies $(\beta)$ or the dual property $(\delta)$
on $\mathbb{C}\setminus S$ and for which there is a bounded open set
$V\subset \mathbb{C}$ with $S\cap V=\emptyset$ or $S\subset V$, then
the following holds:
\begin{enumerate}
\item If $\sigma(T)$ is dominating in $V$ than $T$ has a nontrivial
  invariant subspace.
\item If the essential spectrum of $T$ is dominating in $V$ then the lattice
  $\operatorname{Lat}(T)$ of invariant subspaces is rich. 
\end{enumerate}
Even more generally, the authors showed in \cite{JE-Pru02} a corresponding result
for operators $T$ for which the localizable spectrum of $T$ or of $T^\prime$
is thick (in the above mentioned sense), has a non-trivial invariant subspace.
Notice, that some thickness condition on the spectrum is necessary, as there
exists an example due to Charles Read \cite{Read97} of a quasinilpotent
and hence even decomposable operator on some Banach space without any
non-trivial invariant subspaces. 

Following the success of the Scott Brown technique for single operators,
J\"org set out to extend these ideas to tuples of commuting operators.
The natural goal is to prove that any contractive tuple of commuting operators with sufficiently
rich spectrum has non-trivial invariant subspaces.
In the multivariate setting, several difficulties arise.
As mentioned before, spectral theory becomes substantially more difficult.
Moreover, just like there is more than one generalization of the unit disc to higher variables,
there are several reasonable notions of contractive tuples of operators.
Finally, the Sz.-Nagy--Foia\c{s} $H^\infty$-functional calculus, which is a crucial ingredient
in the classical Scott Brown technique, is generally not available in the multivariable setting.
In \cite{Eschmeier97}, J\"org proved that every commuting row contraction that possesses a spherical dilation
and whose Harte spectrum is dominating in the ball admits non-trivial invariant subspaces.
Along the way, he established an $H^\infty(\mathbb{B}_d)$-functional calculus,
which itself inspired further research approximately 20 years later \cite{CD16a,BHM17}.
He also established a version of this theorem on the polydisc \cite{Eschmeier01}.

In the one variable setting, the Scott Brown technique was refined by Olin and Thomson \cite{OT80}
to show that a subnormal operator $T$ not only admits nontrivial invariant subspace, but is even reflexive.
This means that the weak operator topology closed algebra generated by $T$ can be completely recovered
from its invariant subspace lattice. J\"org successfully established variants
of this result in the multivariable setting \cite{Eschmeier99,DE05,Eschmeier06,Eschmeier01a,Eschmeier07}.
The difficulty of this subject can perhaps be appreciated from the fact that the question of whether every subnormal operator tuple is reflexive remains open to this day.
\begin{figure}
\includegraphics[width=12cm]{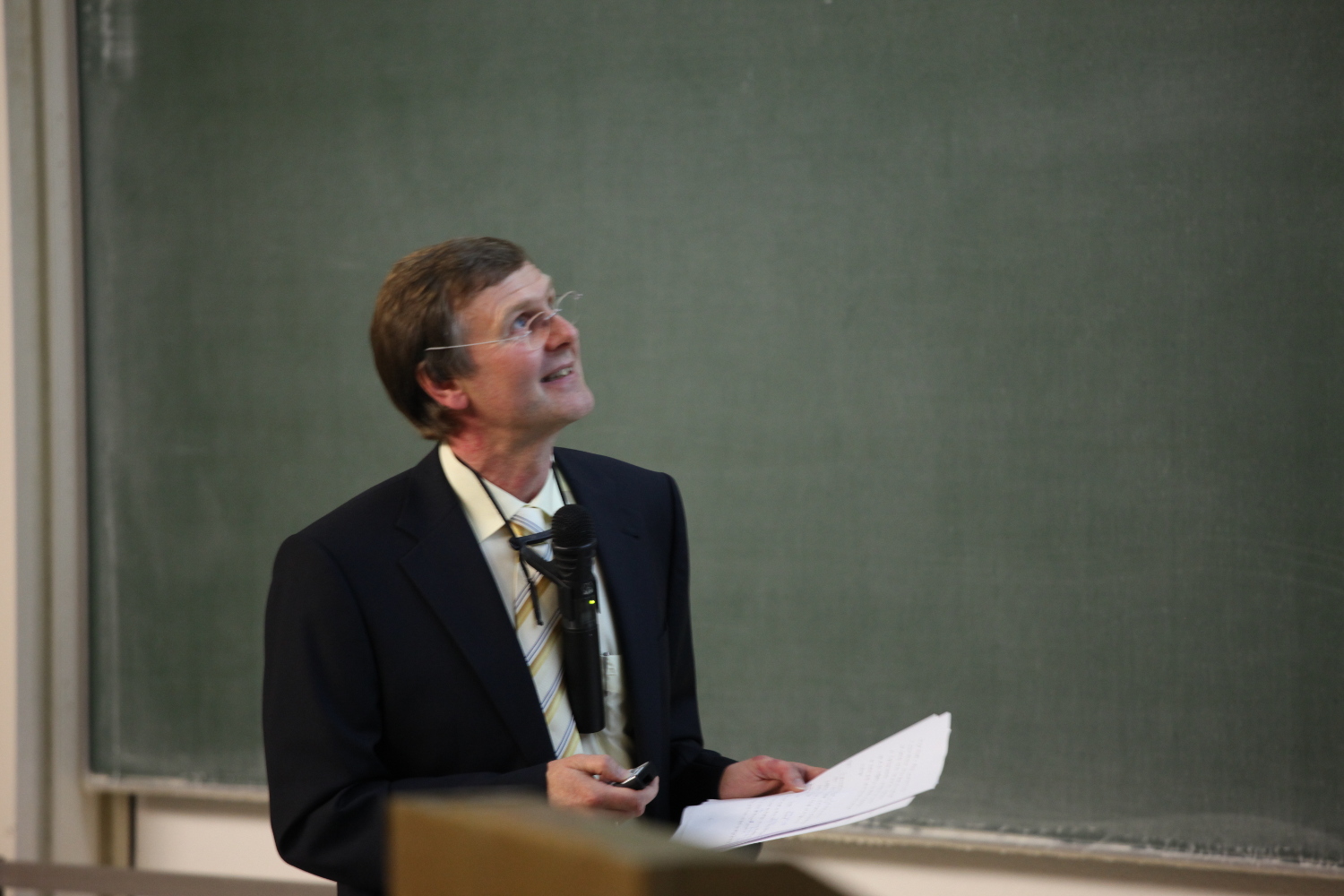}
\caption{J\"org in 2009, at the occasion of EA's 65th birthday.}
\end{figure}

\section{Multivariable operator theory}

The interplay between functional analysis and function theory of several complex variables entered around 1970-ies into a new era thanks to Joseph Taylor's innovative insights,
in particular his novel joint spectrum and the topological-homological approach to functional calculi \cite{Taylor-1970,Taylor-1972}. J\"org was a central figure and inspired contributor to this new chapter of modern analysis. His many notable contributions would fill an entire volume, way beyond the size of the present biographical note. We outline a couple of snapshots.

In duality with the fibre product of complex spaces, the tensor product of analytic modules is inviting at analyzing the blend of spectral behaviors of its factors. Two traditional constructions of operator theory stand aside in this respect. More specifically, consider a pair of Banach spaces $X,Y$ and commutative tuples of linear bounded operators acting on them
$S \in L(X)^m$ and $T \in L(Y)^n$. With a fixed cross-norm, the completed tensor product $X \otimes Y$ carries the commutative $(n+m)$-tuple $\tau = (S \otimes I, I \otimes T)$.
One of J\"org's early works \cite{JE-1988} (see also the
  third chapter in his Habilitation Thesis \cite{JE-1987})
offers a complete evaluation of the joint spectrum of $\tau$, its joint essential spectrum, and the values of the Fredholm index.
A second natural path of melting tuples of linear transforms into new ones has to do with the so-called elementary operators. Namely, let $J$ be a bilateral ideal in the space of linear bounded operators $L(K,H)$ acting between two Hilbert spaces $K,H$. Let $S \in L(H)^n$ and $T \in L(K)^n$ be tuples of commuting elements. The operator
$$ R : J \longrightarrow J, \ \ R(A) = \sum_{j=1}^n S_j A T_j,$$
is an elementary operator.  A complete spectral picture of $R$ has quite surprising implications, for instance to solving non-commutative operator equations relevant to control and stability theory of systems. The spectrum, essential spectrum, and Fredholm index of $R$ were first computed by one of us (R.E.C.) together with Lawrence A. Fialkow \cite{CuFi}, and further refined under the assumption of an intricate geometric condition, known as the ``finite fibre property" \cite{Fialkow-1986}. The note \cite{JE-MP-1995} shows, employing cohomology and analytic localization techniques, that the finite fibre property is always satisfied. The investigations of tensor products of analytic modules and separately of elementary operators continue to flourish.

As defined, Taylor's joint spectrum of a commutative tuple of linear operators is quite elusive, requiring a grasp of higher torsion spaces of a pair of analytic modules. Cohomologically trivial cases are in this sense a treat, as much as, keeping the proportion, pseudoconvex domains are among all domains in $\C^d$.
For this reason, a collection of simple examples is both precious and inspiring. Take for instance a bounded pseudoconvex domain $\Omega$ in $\C^d$ and a finite system of bounded analytic functions defined on $\Omega$: $f = (f_1, f_2, \ldots, f_n)$. The note \cite{JE-MP-1993} contains a description of Taylor's joint spectrum of $f$, as acting on Bergman space $L_a^2(\Omega)$: 
$$
 \sigma(f, L^2_a(\Omega)) = \overline{f(\Omega)}.
$$
In the same setting, assuming the boundary of $\Omega$ is smooth, strictly pseudoconvex, the joint essential spectrum is
$$ 
\sigma_e(f, L^2_a(\Omega)) = \bigcap_U \overline{ f(U \cap \Omega)},$$
where $U$ runs over all open neighborhoods of $\partial \Omega$. These results can be interpreted as solutions to division problem with $L^2$-bounds, in the spirit of the celebrated Corona Problem.

Bounded analytic interpolation in several complex variables is significantly more challenging than the classical Nevanlinna-Pick or Carath\'eodory-Fej\'er 1D problems. A major breakthrough was recorded in the 1990-ies with the isolation of the Schur-Agler class of analytic multipliers on a Hilbert space of analytic functions with a reproducing kernel. It was Jim Agler who recognized at that time that this new algebra provides the correct analog for stating and proving multivariate analogs of the quasi-totality of known bounded analytic interpolation results in one dimension \cite{Agler-1990}. J\"org entered into the first line of avant-garde researchers charting this new territory, as for instance, the articles \cite{JE-LP-MP,JE-MP-2002,JE-MP-2003} amply illustrate. 

A commuting tuple $T = (T_1,\ldots,T_n)$ of operators is said to be Fredholm if the cohomology groups $H^p(T)$ of the Koszul complex are
all finite dimensional. The dimensions of $H^p(T)$ and related objects naturally carry operator theoretic meaning.
It turns out that they can also be related to certain analytical quantities.
In his paper \cite{JE-2007-bis}, J\"org showed that the limits $\lim_{k \to \infty} \dim H^p(T^k)/k^n$
exist and in fact agree with the so-called Samuel multiplicities of stalks of cohomology sheaves;
moreover, they agree with the generic dimension of $H^p(z- T)$ near $z=0$.
This result employs a beautiful blend of operator theory, commutative algebra and analytic geometry.
J\"org continued his investigations in \cite{JE-2007,Eschmeier08}.

The model theory of Sz.-Nagy and Foia\c{s} yields a complete unitary invariant for (completely non-unitary) contractions,
namely the characteristic function. This is closely related to the theory of the Hardy space on the disc and classical inner functions. J\"org, together with several coauthors, made advances into a corresponding multivariable theory \cite{BES06} and into setting of more general reproducing kernel Hilbert spaces \cite{Eschmeier18,ET21,BEK+17}.

Another topic in which J\"org was active is the characterization of the essential commutant of the analytic Toeplitz operators. A theorem of Davidson \cite{Davidson77} shows that, on the Hardy space, an operator $T$ commutes modulo the compacts with every analytic Toeplitz operator if and only if $T$ is a compact perturbation of a Toeplitz operator with symbol in $H^\infty + C$. Very general extensions to the multivariable setting were given by Eschmeier and coauthors in \cite{EE15,DEE11,D-JE-S-2017}.
\begin{figure}
\includegraphics[width=12cm]{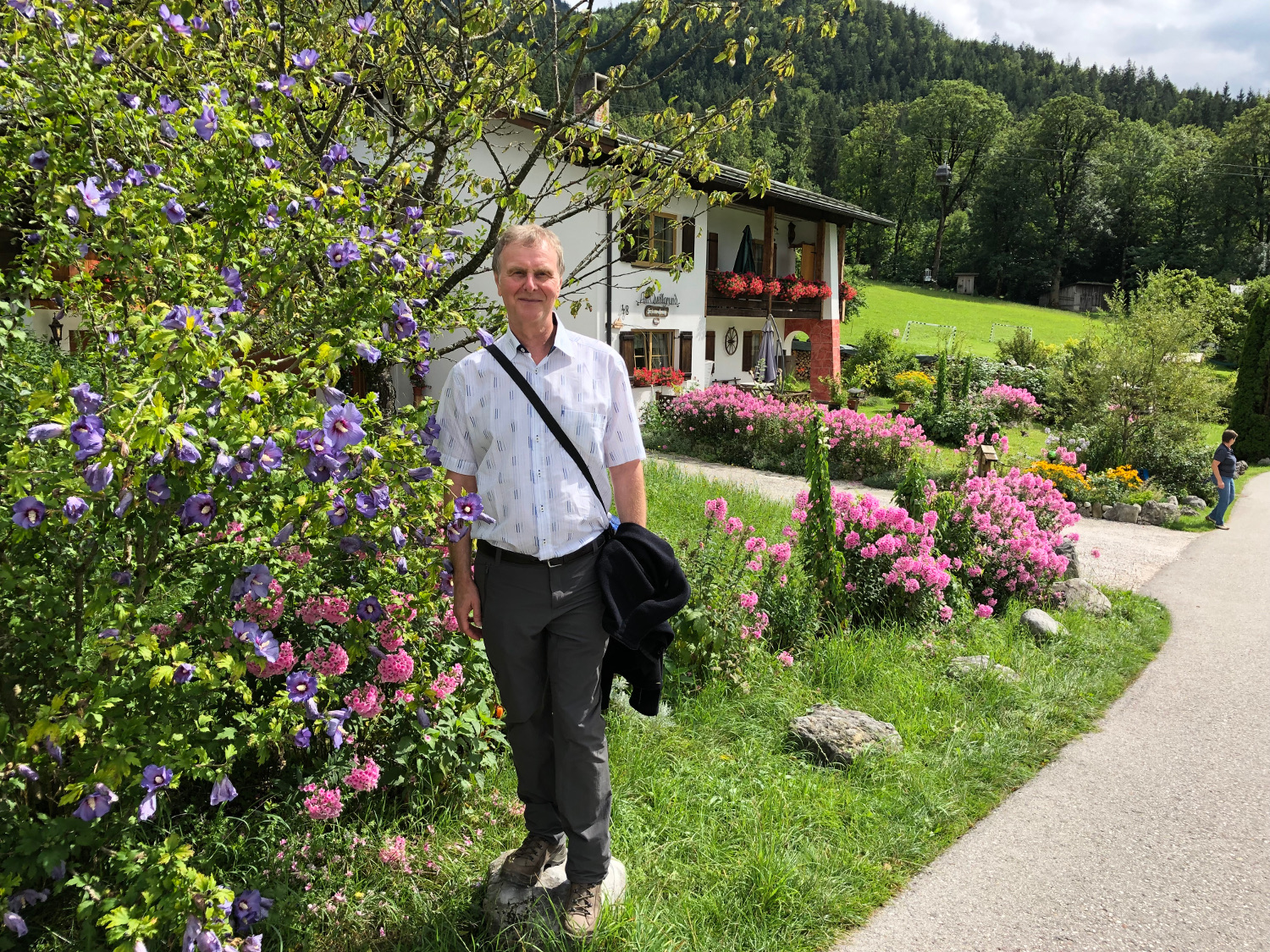}
\caption{J\"org in 2019}
\end{figure}

\section{Arveson-Douglas conjecture}

Classically, there is an intimate relationship between the study of contraction operators
on Hilbert space and the Hardy space $H^2$ on the unit disc.
In the theory of tuples of commuting operators on a Hilbert space,
the appropriate generalization of the Hardy space is the Drury--Arveson space $H^2_d$,
see \cite{Arveson98}.
This is the reproducing kernel Hilberts space of analytic functions on the
Euclidean unit ball $\mathbb{B}_d$ in $\mathbb{C}^d$ with reproducing kernel
\begin{equation*}
  K(z,w) = \frac{1}{1 - \langle z,w \rangle }.
\end{equation*}
Given J\"org's experience in multivariable operator theory and in several complex variables,
it is no surprise that he made important contributions to this subject early on \cite{JE-MP-2002}.

Since the coordinate functions $z_1,\ldots,z_d$ are multipliers of $H^2_d$,
the space $H^2_d$ becomes a module over the polynomial ring.
To each (closed) submodule $M \subset H^2_d$, one associates the operators
\begin{equation*}
  S_j: M^\bot \to M^\bot, \quad f \mapsto P_{M^\bot} (z_j f),
\end{equation*}
which are commuting linear operators on $M^\bot$.
A very influential conjecture, made by Arveson \cite{Arveson02a} and refined by Douglas \cite{Douglas06},
asserts that if $M$ is the closure of a homogeneous ideal $I$ of polynomials, then the cross commutators
$[S_j , S_k^*]$ belong to the Schatten class $\mathcal{S}^p$ for all $p > \dim Z(I)$, where $Z(I)$
denotes the zero set of $I$.
This conjecture has attracted a large amount of attention.

The initial motivation for the Arveson--Douglas conjecture came from Arveson's work on the curvature invariant
of certain operator tuples \cite{Arveson02a}, but the conjecture turned out to be interesting for other reasons as well.
For instance, if $I$ is a homogeneous ideal of infinite co-dimension
so that the Arveson--Douglas conjecture holds for $M = \overline{I}$,
then the quotient of the Toeplitz $C^*$-algebra associated with $M$ modulo the compacts is commutative.
In fact, the quotient is isomorphic to $C(Z(I) \cap \partial \mathbb{B}_d)$.
As Douglas observed, this gives rise to a $K$-homology element of the space $Z(I) \cap \partial \mathbb{B}_d$ \cite{Douglas06}.

After the Arveson--Douglas conjecture was formulated, it was verified in a number of special cases.
J\"org reduced the conjecture to a certain operator inequality \cite{JE-2011}.
This led to a unified proof of all cases in which the conjecture was known to hold at the time.
A few years later, J\"org Eschmeier and Miroslav Engli\v{s} achieved a spectacular breakthrough
by showing that the conjecture holds if the homogeneous ideal $I$ is the vanishing ideal
of a homogeneous variety that is smooth away from the origin \cite{E-JE-2015}.
A similar, related result was independently obtained by Douglas, Tang and Yu around the same time \cite{DTY16}.
These works have inspired a lot of research, and the area remains very active to this day.

\section{Teaching and mentoring}

J\"org Eschmeier was a very highly regarded teacher and mentor.
His lectures were widely known for their clarity and precision.
They were frequently attended by a large number of students, and even
his more advanced courses in functional analysis and complex analysis drew in students from outside of mathematics.
At Saarland University, he won the award for the best lecture in mathematics
five times, more often than any of his colleagues in the
  department.

In addition, J\"org was an extremely dedicated mentor for bachelor's, master's and doctoral theses.
He was very generous with his time and invested considerable effort into advising students.

A particular gem among J\"org's lectures were those about several complex variables.
His lecture notes formed the basis for his book \cite{JE-book}.
The book stands out in that it provides a self-contained
treatment of important theorems in several complex variables,
assuming only basic undergraduate analysis as a prerequisite.
The book starts with basic properties of holomorphic functions in several variables, treats analytic sets
and domains of holomorphy, proves Oka's theorems and explains the solution of the Levi problem.
The last chapter contains beautiful applications to functional analysis, namely the Arens-Calder\'on functional calculus,
Shilov's idempotent theorem and the Arens-Royden theorem.
Through his lectures and his book, J\"org made this fascinating subject accessible to many students.

Another very nice example of informative and well structured writing are his lectures on invariant
  subspaces contained in Part III of \cite{ADEL}.

Already during his period as an Alexander von Humboldt Fellow, J\"{o}rg
was involved in the mentoring of Roland Wolff, a doctoral student of George
Maltese, whom he introduced to the field of Bergman and Hardy spaces in several
variables. Roland Wolff finished his thesis
\textit{Spectral theory on Hardy spaces in several complex variables}
\cite{Wolff-Th}
in Saarbr\"ucken, where he became the first assistent of J\"org. Part of this
theses is also published in \cite{Wolff-1}.

The following is a list of all completed doctoral theses written under the
mentorship of J\"org at Saarbr\"ucken:
\begin{enumerate}
\item Michael Didas, \textit{On the structure of von Neumann $n$-tuples
    over strictly pseudoconvex sets} \cite{Didas-Th}. 
  Michael Didas considers operator tuples $T=(T_1,\ldots,T_n)$ on a Hilbert
  space $H$ admitting a contractive functional calculus
  $\Phi_T:A(D)\to L(H)$, where $D$ is a strictly pseudoconvex subset of a
  Stein submanifold of $\mathbb{C}^n$ which admit a $\partial D$-unitary
  dilation.  Using dual algebra methods he obtains
  reflexivity results and invariant subspace results for
  large classes commuting 
  operator tuples $T=(T_1,\ldots,T_n)$ on a Hilbert space with dominating
  Harte- or Taylor spectra in $D$. 
  The main parts of this thesis have been published in \cite{Didas-1}.
  Though
  Michael Didas left the university he always stayed in contact with J\"{o}rg
  and they published seven joint articles. 
\item Eric R\'{e}olon, \textit{Zur Spektraltheorie vertauschender Operatortupel:
    Fredholmtheorie und subnormale Operatortupel} \cite{Reolon-Th}.
  In his thesis Eric R\'{e}olon obtains a Banach space variant of a
  Fredholm index formula for essentially normal tuples given in
  \cite{JE-MP-book}, Theorem 10.3.15. He also shows that an operator
  tuple on a Hilbert space is essentially subnormal if and only if it has an
  essentially normal extension and if and only if it has an extension to a
  compact perturbation of a normal tuple. For single operators this had been
  shown by N.S.~Feldman in \cite{Fe99}.
\item Christoph Barbian, \textit{Beurling-Type Representation of
    Invariant Subspaces in Reproducing Kernel Hilbert Spaces}
  \cite{Barbian-Th}.
  Christoph Barbian studies invariant subspaces of reproducing kernel Hilbert spaces, including
  the difficult case of the Bergman space on the unit disc.
  He introduces the notion of Beurling decomposability, and obtains
  criteria for this property to hold.
  Among other things, this has implications for multivariable spectra of multiplication tuples.
  This thesis as well as the following ones is available
  for download.
  For further developments see also
  \cite{Bar08,Bar09,Bar11}. 
\item Dominik Faas, \textit{Zur Darstellungs- und Spektraltheorie für
    nichtvertauschende Operatortupel} \cite{Faas-Th}. This thesis is related
  to J\"{o}rg's work on Samuel multiplicities. A local closed range theorem
  for semi-Fredholm valued functions was later improved to a global version
  by Dominik Faas and J\"{o}rg in \cite{EF10}
\item Kevin Everard, \textit{A Toeplitz projection for multivariable isometries}
  \cite{Everard-Th}. For compact sets $K\subset \mathbb{C}^n$ and closed
subalgebras $A$ of $C(K)$ J\"{o}rg introduced in \cite{Eschmeier06} the notion
of an $A$-isometry. This class of commuting $n$-tuples of operators
includes spherical isometries and $n$-tuples of commuting isometries.
This allowed J\"{o}rg together with Michael Didas and Kevin Everard to
introduce associated analytic Toeplitz operators \cite{DEE11}. The thesis
of Kevin Everard completes that approach (see also \cite{EE15}).
\item Michael Wernet, \textit{On semi-Fredholm theory and essential normality}
  \cite{Wernet-Th}. Michael Wernet contributes to four areas of J\"{o}rg's
  interests: He extends J\"{o}rg's results of \cite{Esc08a}, he shows that a
  number of positive results on the Arveson-Douglas conjecture can be extended
  to arbitrary graded Hilbert modules and that the validity of the conjecture
  is equivalent for a large class of analytic functional Hilbert spaces, he
  generalizes an essential von Neumann inequality of Matthew Kennedy and Orr
  Shalit (\cite{KS15}, Theorem 6.1), and using a result of \cite{EF10}
  answers a question by Ronald Douglas (\cite{Dou09}, Question 1)
  and (extending some results of J\"{o}rg and Johannes Schmitt \cite{ES14})
  gives a partial answer to another question of Ronald Douglas in
  \cite{Dou09}, Question 3.
\item Dominik Schillo, \textit{K-contractions, and perturbations of Toeplitz
    operators} \cite{Schillo-Th}. For many analytic Hilbert
  function spaces with reproducing kernels $K$ Dominik Schillo obtains model
  theorems for $K$-contractions. In a second part he studies Toeplitz
  operators associated with regular $A$-isometries and uses
  methods from \cite{DEE11} to characterize finite rank and Schatten-p-class
  perturbations of analytic Toeplitz operators. Part of these results have
  also been published in \cite{D-JE-S-2017}. 
\item Sebastian Langend\"orfer, \textit{On unitarily invariant spaces and
    Cowen-Douglas theory} \cite{Langend-Th}. A rather general version of a
  characterization of Toeplitz operators with pluriharmonic symbols on
  unitarily invariant with appropriate reproducing kernel and an extension
  of results by Chang, Chen and Fang \cite{CCF}
  to the several variable case are given. See also the joint works
  \cite{EL17,EL18,EL19} of Sebastian Langend\"{o}rfer with J\"{o}rg.
\item Daniel Kraemer, \textit{Toeplitz operators on Hardy spaces}
  \cite{Kraemer-Th}. A several variable Toeplitz operator theory is
  developed on Hardy type $H^p(G)$ spaces which is applicable for
  bounded symmetric domains and bounded strictly pseudoconvex domains.
  In particular, in the situation of strictly pseudoconvex domains he
  obtains a generalization of J\"{o}rg's spectral mapping theorem
  from \cite{Esc13b}.
\end{enumerate}

The high quality of these theses reflects the outstandig quality of J\"{o}rg
Eschmeier as an academic teacher.\medskip

\subsection*{Acknowledgement} We are grateful to Steliana Eschmeier for her support in the preparation of this article.
  \nocite{*}

  \printbibliography[title=Complete list of J\"org Eschmeier's works,keyword={JEarticle}]
  \printbibliography[title=Other References,notkeyword=JEarticle]

\end{document}